\documentstyle[amscd, amstex]{amsart}

\def\NN        {{\Bbb N}}
\def\ZZ         {{\Bbb Z}}
\def\RR         {{\Bbb R}}
\def\CC         {{\Bbb C}}
\def\QQ         {{\Bbb Q}}

\newtheorem{thm}{Theorem}[section]
\newtheorem{lem}[thm]{Lemma}
\newtheorem{cor}[thm]{Corollary}
\newtheorem{pr}[thm]{Proposition}
\newtheorem{conj}[thm]{Conjecture}
\theoremstyle{definition}
\newtheorem{rem}[thm]{Remark}

\newtheorem{defn}[thm]{Definition}

\newcommand{\cv}{{\cal V}}
\newcommand{\pro}{{\rm Proj}}

\newcommand{\xd}{{{ X}_{\Delta}}}
\newcommand{\sd}{{ \Sigma_\Delta}}
\newcommand{\xsd}{{{X}_{ \Sigma_\Delta}}}

\newcommand{\co}{{\rm Conv}}

\newcommand{\spe}{{\rm Spec}}
\newcommand{\key}{\bibitem}

\newcommand{\xs}{{{ X}_{\Sigma}}}

\newcommand\hidot{{\raise1pt\hbox{$\scriptscriptstyle\bullet$}}}
\newcommand\lodot{{\raise.3pt\hbox{$\scriptscriptstyle\bullet$}}}

\newcommand{\inte}{{\rm int}}

\newcommand{\Hom}{{\rm Hom}}

\newcommand{\ve}{{\scriptscriptstyle\vee}}
\begin{document}

\title[  Mirror Symmetry  for Calabi-Yau   complete intersections]
{ Mirror Symmetry  for Calabi-Yau \\   complete intersections in
Fano toric varieties}
\author{Anvar R. Mavlyutov}
\address {Department of Mathematics, Oklahoma State  University, Stillwater, OK, USA.}
 \email{mavlyutov@@math.okstate.edu}



\maketitle

\begin{abstract}
Generalizing the notions of reflexive polytopes and
nef-partitions of  Batyrev and Borisov, we propose a mirror
symmetry construction for Calabi-Yau complete intersections in
Fano toric varieties.
\end{abstract}

\tableofcontents

 \setcounter{section}{-1}

 \section{Introduction.}

  Reflexive polytopes $\Delta$ in $\RR^d$,
introduced by Victor Batyrev in \cite{b}, are determined by the
property that they have vertices in $\ZZ^d$ and have the origin in
their interior
 with the polar dual  polytope $$\Delta^*=\{y\in  \RR^d\mid \langle \Delta,y\rangle\ge-1\}$$    satisfying
the same property. The polar duality gives an involution between
the sets of reflexive polytopes: $(\Delta^*)^*=\Delta$. There is a
one-to-one correspondence between isomorphism classes of reflexive
polytopes $\Delta$ in $\RR^d$ and $d$-dimensional Gorenstein Fano
toric varieties  given by $$\Delta\mapsto
X_\Delta:=\pro(\CC[\RR_{\ge0}(\Delta,1)\cap\ZZ^{d+1}]),$$  where
the grading is induced by the last coordinate in $\ZZ^{d+1}$. The
dual pair of reflexive polytopes $\Delta$ and $\Delta^*$
corresponds to the {\it Batyrev mirror pair} of ample Calabi-Yau
hypersurfaces $Y_\Delta\subset X_\Delta$ and $Y_{\Delta^*}\subset
X_{\Delta^*}$ in Gorenstein Fano toric varieties in \cite{b}. By
taking maximal projective  crepant partial resolutions $\widehat
Y_\Delta\rightarrow Y_\Delta$ and $\widehat
Y_{\Delta^*}\rightarrow Y_{\Delta^*}$ induced by toric blow ups,
Batyrev obtained a mirror pair of minimal Calabi-Yau hypersurfaces
$\widehat Y_\Delta,\widehat Y_{\Delta^*}$.

Generalizing the polar duality of reflexive polytopes, Lev Borisov
in \cite{bo} introduced the notion of {\it nef-partition}, which
is a Minkowski sum decomposition
 of the reflexive polytope
 $\Delta=\Delta_1 +\cdots+\Delta_r$ by lattice polytopes
 such that the origin $0\in\Delta_i$ for all $1\le i\le r$.
A nef-partition has a {\it dual nef-partition} defined as the
Minkowski sum decomposition of the reflexive polytope
$\nabla=\nabla_1+\dots+\nabla_r$ in the dual vector space with
 $\nabla_j$   determined by  $\langle
\Delta_i,\nabla_j\rangle\ge-\delta_{ij}$ for all $1\le i,j\le r$,
where $\delta_{ij}$ is the Kronecker symbol.

A nef-partition of a reflexive polytope
$\Delta=\Delta_1+\cdots+\Delta_r$ with $r<d$ and $\dim\Delta_i>0$,
for all $1\le i\le r$, defines a nef Calabi-Yau complete
intersection $Y_{\Delta_1,\dots,\Delta_r}$ in the  Gorenstein Fano
toric variety $X_\Delta$ given by the equations:
$$
 \Biggl(\sum_{m\in  \Delta_i\cap \ZZ^d }a_{i,m}
\prod_{v_\rho\in\cv(\Delta^*)}x_{\rho}^{ \langle m,v_\rho\rangle
}\Biggr)\prod_{
 v_\rho \in  \nabla_i  } x_\rho=0,\quad i=1,\dots,r,$$
with generic $a_{i,m}\in\CC$, where $x_\rho$ are the Cox
homogeneous coordinates of the toric variety $X_\Delta$
corresponding   to the vertices $v_\rho$ of the polytope
$\Delta^*$.

The Batyrev-Borisov mirror symmetry construction is a pair of
families of generic nef Calabi-Yau complete intersections
$Y_{\Delta_1,\dots,\Delta_k}\subset X_\Delta$ and
$Y_{\nabla_1,\dots,\nabla_k}\subset X_{\nabla}$ in Gorenstein Fano
toric varieties corresponding to a dual pair of nef-partitions
$\Delta=\Delta_1+\cdots+\Delta_r$ and
$\nabla=\nabla_1+\cdots+\nabla_r$. By taking maximal projective
crepant partial resolutions $\widehat
Y_{\Delta_1,\dots,\Delta_k}\rightarrow
Y_{\Delta_1,\dots,\Delta_k}$ and $\widehat
Y_{\nabla_1,\dots,\nabla_k} \rightarrow
Y_{\nabla_1,\dots,\nabla_k}$, one obtains the {\it Batyrev-Borisov
mirror pair} of minimal Calabi-Yau complete intersections.

The topological mirror symmetry test for compact $n$-dimensional
Calabi-Yau manifolds $V$ and $V^*$ is a symmetry of their Hodge
numbers: $$h^{p,q}(V)=h^{n-p,q}(V^*), \quad 0\le p,q\le n.$$ For
singular varieties Hodge numbers must be replaced by the stringy
Hodge numbers $h^{p,q}_{\rm st}$ introduced by V. Batyrev in
\cite{b2}. The usual Hodge numbers coincide with the stringy Hodge
numbers for nonsingular Calabi-Yau varieties, and all crepant
partial resolutions $\widehat V$ of singular Calabi-Yau varieties
$V$ have the same stringy Hodge numbers: $h^{p,q}_{\rm
st}(\widehat V)=h^{p,q}_{\rm st}(  V )$. In \cite{bb2}, Batyrev
and Borisov show that the pair of generic Calabi-Yau complete
intersections $V= Y_{\Delta_1,\dots,\Delta_r}$  and
$V^*=Y_{\nabla_1,\dots,\nabla_r}$ pass the topological mirror
symmetry test:
$$h^{p,q}_{\rm st}(Y_{\Delta_1,\dots,\Delta_r})=h^{d-r-p,q}_{\rm
st}(Y_{\nabla_1,\dots,\nabla_r}),\quad 0\le p,q\le d-r.$$

Generalizing the notions of reflexive polytopes and nef-partitions
for rational polytopes we introduce the notions of $\QQ$-reflexive
polytopes and $\QQ$-nef-partitions. A {\it $\QQ$-reflexive
polytope} $\Delta$  in $\RR^d$ is determined by the properties
    that    $0$ is in the interior of $\Delta$
and $$\co( (\co(\Delta\cap \ZZ^d))^*\cap \ZZ^d)=\Delta^*.$$
 We
show that a $\QQ$-reflexive polytope $\Delta$ corresponds to the
Fano toric variety $X_\Delta$ with at worst canonical
singularities. A $\QQ$-reflexive polytope has a {\it dual
$\QQ$-reflexive polytope} defined by
$\Delta^\circ:=(\co(\Delta\cap \ZZ^d))^*$ and the property
$(\Delta^\circ)^\circ=\Delta$ gives an involution on the set of
$\QQ$-reflexive polytopes.
 The dual lattice polytope $\Delta^*$ of a
$\QQ$-reflexive polytope is called an {\it almost reflexive
polytope} and there is a similar involution on the set  of almost
reflexive polytopes. All reflexive polytopes are $\QQ$-reflexive
and almost reflexive.

 A {\it
$\QQ$-nef-partition} is a Minkowski sum decomposition
$\Delta=\Delta_1+\cdots+\Delta_r$ of a $\QQ$-reflexive polytope
into polytopes in $\RR^d$ such that
  $0\in\Delta_i$ for all $1\le i\le r$, and $$\co(\Delta\cap \ZZ^d)=\co(\Delta_1\cap \ZZ^d)+\cdots+\co(\Delta_r\cap
  \ZZ^d).$$
We prove that  a $\QQ$-nef-partition
$\Delta=\Delta_1+\cdots+\Delta_r$ has a {\it dual
$\QQ$-nef-partition} $\nabla=\nabla_1+\cdots+\nabla_r$ determined
by $\langle \Delta_i\cap \ZZ^d,\nabla_j\rangle\ge-\delta_{ij}$ for
all $1\le i,j\le r$.  This  dual pair of   $\QQ$-nef-partitions
corresponds to a pair of $\QQ$-nef Calabi-Yau complete
intersections $Y_{\Delta_1,\dots,\Delta_k}\subset X_\Delta$ and
$Y_{\nabla_1,\dots,\nabla_k}\subset X_{\nabla}$ in Fano toric
varieties. We expect that this pair passes the topological mirror
symmetry test as in \cite{bb2}.

In \cite{bb}, Batyrev and Borisov introduce the notion of
reflexive Gorenstein cones $\sigma\subset\RR^{\bar d}$, which
canonically correspond to  Gorenstein Fano toric varieties
$X_\sigma=\pro(\CC[\sigma^\ve\cap \ZZ^{\bar d}])$ such that ${\cal
O}_{X_\sigma}(1)$ is an ample invertible sheaf and there is a
positive integer $r$ such that ${\cal O}_{X_\sigma}(r)$ is
isomorphic to  the anticanonical sheaf of $X_\sigma$. The zeros
$Y_\sigma$ of generic global sections of ${\cal O}_{X_\sigma}(1)$
are called {\it generalized Calabi-Yau manifolds}.
 The dual cone $$\sigma^\ve =\{y\in \RR^{\bar d}\mid \langle x, y\rangle\ge0
\,\forall\,  x\in\sigma \}$$  of a reflexive Gorenstein cone
$\sigma$ is again reflexive, and the dual pair $\sigma$ and
$\sigma^\ve$ corresponds to the {\it mirror pair} of generalized
Calabi-Yau manifolds $Y_\sigma$ and  $Y_{\sigma^\ve}$, which are
ample   hypersurfaces in the respective Gorenstein Fano toric
varieties.

Combining the ideas of \cite{bb} with the notion of almost
reflexive polytopes, we introduce the notion  of {\it almost
reflexive} Gorenstein cones $\sigma$. Their dual cones
$\sigma^\ve$ are no longer Gorenstein, but there is a canonically
defined grading on $\sigma^\ve\cap \ZZ^{\bar d}$. This allows us
to associate to an almost reflexive Gorenstein cone
$\sigma\subset\RR^{\bar d}$ the Fano toric variety
$X_\sigma=\pro(\CC[\sigma^\ve\cap \ZZ^{\bar d}])$. The reflexive
rank one sheaf ${\cal O}_{X_\sigma}(1)$ corresponds to an ample
$\QQ$-Cartier divisor and there is a positive integer $r$ such
that ${\cal O}_{X_\sigma}(r)$ is isomorphic to the anticanonical
sheaf on $X_\sigma$. In particular, we have  a   generalized
Calabi-Yau manifold $Y_\sigma$  given by generic global sections
of ${\cal O}_{X_\sigma}(1)$.  There is an involution on the set of
almost reflexive Gorenstein cones $\sigma\mapsto\sigma^{\bullet}$.
For a dual pair of  almost reflexive Gorenstein cones $\sigma$ and
$\sigma^\bullet$ we expect that the correspondence  between
generalized Calabi-Yau manifolds $Y_\sigma \leftrightarrow
Y_{\sigma^\bullet}$ corresponds to the mirror involution in $N=2$
super conformal field theory.

\section{$\QQ$-reflexive and almost reflexive polytopes.} \label{s:refl}

In this section, we first review the definition   of reflexive
polytopes  due to V. Batyrev in \cite{b}, and then construct a
natural generalization of these notions for rational and lattice
polytopes.

 Let $M$ be a lattice of rank $d$ and $N=\Hom(M,\ZZ)$ be its dual lattice with
  a natural paring $\langle \underline{\hskip0.2cm},\underline{\hskip0.2cm}\,\rangle:M\times
  N\rightarrow\ZZ$. Denote $M_\RR=M\otimes_\ZZ \RR$,
  $N_\RR=N\otimes_\ZZ \RR$.

\begin{defn}
  A $d$-dimensional lattice polytope $\Delta\subset M_\RR$  is called a {\it canonical Fano polytope}
 if ${\rm int} (\Delta)\cap M=\{0\}$.
\end{defn}

\begin{defn}\cite{b}\label{d:refbat} A $d$-dimensional lattice polytope $\Delta\subset M_\RR$
    is called  {\it reflexive} (with respect to $M$) if
$0\in {\rm int} (\Delta)$ and   the  dual polytope
$$\Delta^*=\{n\in N_\RR\mid \langle m,n\rangle\ge-1\, \forall\, m\in \Delta\}$$
in the dual vector space $N_\RR$ is also a lattice polytope. The
pair $\Delta$ and $\Delta^*$ is called a  {\it dual pair reflexive
polytopes} and it satisfies $\Delta=(\Delta^*)^*$.
\end{defn}

\begin{defn} A  compact toric  variety $X$ is called\\
  $\bullet$ {\it Fano} if the anticanonical divisor  $-K_X$ is
ample and $\QQ$-Cartier, \\ $\bullet$ {\it Gorenstein} if $K_X$ is
Cartier.
\end{defn}

\begin{pr} There is a bijection between isomorphism classes
of canonical Fano polytopes and    Fano toric varieties with
canonical singularities given by $\Delta\mapsto X_{\Delta^*}$. In
particular, Gorenstein Fano toric varieties correspond to
reflexive polytopes.
\end{pr}

Generalizing the notion  of a reflexive polytope we introduce:

\begin{defn} A $d$-dimensional  polytope $\Delta$ in $M_\RR$
    is called  {\it $\QQ$-reflexive} (with respect to $M$) if
   $0\in {\rm int} ( \Delta )$
and
\begin{equation}\label{e:qreflex}\co( (\co(\Delta\cap M ))^*\cap N)=\Delta^*.\end{equation}
\end{defn}

\begin{rem} \label{r:refqref} For a $\QQ$-reflexive polytope $\Delta$ its dual   $\Delta^*$  is a lattice polytope,
whence reflexive polytopes are  $\QQ$-reflexive. It follows from
$(\ref{e:qreflex})$ that a $\QQ$-reflexive polytope is rational,
i.e., its vertices lie in $M_\QQ$. These properties together with
the next ones suggest  the name of $\QQ$-reflexive.
\end{rem}

\begin{defn}\label{d:circ} Denote $[\Delta]:=\co(\Delta\cap M)$ for a
polytope $\Delta$ in $M_\RR$ (and, similarly in $N_\RR$). Also,
define $\Delta^\circ:= [\Delta]^*=(\co(\Delta\cap M))^*$
\end{defn}

 In this
notation, equation $(\ref{e:qreflex})$ is $[[\Delta]^*]=\Delta^*$,
or, equivalently,  $(\Delta^\circ)^\circ=\Delta$. Hence, we have

\begin{lem} If $\Delta\subset M_\RR$
    is  $\QQ$-reflexive , then $\Delta^\circ=(\co(\Delta \cap M))^*\subset N_\RR$ is
    $\QQ$-reflexive and the map $\Delta\mapsto \Delta^\circ$ is  an involution on the set  of $\QQ$-reflexive
    polytopes.
\end{lem}

  We
will call the pair of rational polytopes  $\Delta$ and
$\Delta^\circ$  as {\it the dual pair  of $\QQ$-reflexive
polytopes}.

\begin{rem} A $\QQ$-reflexive polytope $\Delta$ is completely determined
by the convex hull $[\Delta]$ of its lattice points since
$\Delta=[[\Delta]^*]^*$.
\end{rem}

\begin{defn} A $d$-dimensional  lattice polytope $\Delta$ in $N_\RR$
    is called  {\it almost reflexive} (with respect to $N$) if
   $0\in {\rm int} ( \Delta )$
and
\begin{equation}\label{e:almreflex}\co( \co(\Delta^*\cap M))^*\cap N)=\Delta.\end{equation}
\end{defn}

\begin{lem}\label{qrefaref} A   polytope $\Delta$ in $M_\RR$ is $\QQ$-reflexive
if and only if $\Delta^*$ in $N_\RR$ is almost reflexive. In
particular, reflexive polytopes are almost reflexive.
\end{lem}

\begin{defn}\label{d:bul}  For a
polytope $\Delta$ in $N_\RR$  define $\Delta^\bullet:=
[\Delta^*]=\co(\Delta^*\cap M)$
\end{defn}

\begin{rem} \label{r:almref}  In the new notation, equation (\ref{e:almreflex})
is $[[\Delta^*]^*]=\Delta$, or, equivalently,
$(\Delta^\bullet)^\bullet=\Delta$.
\end{rem}

\begin{lem} If $\Delta\subset N_\RR$
    is  almost reflexive , then $\Delta^\bullet:=\co(\Delta^*\cap M)$  is
     almost reflexive and the map $\Delta\mapsto \Delta^\bullet$ is  an involution on the set  of almost reflexive
    polytopes.
\end{lem}

We will call the pair of lattice polytopes  $\Delta$ and
$\Delta^\bullet$  as {\it the dual pair  of  almost reflexive
polytopes}.

A $\QQ$-reflexive polytope  has the following properties.

\begin{lem}\label{l:facet} Every facet of a  $\QQ$-reflexive polytope contains a
lattice point.
\end{lem}

\begin{pf} Suppose that $\Delta$ is $\QQ$-reflexive. Since
$\Delta^*$ is a lattice polytope, every facet  of $\Delta$ is
determined by $\{m\in M_\RR\mid \langle m, v\rangle=-1\}$ for a
vertex $v\in\Delta^*$. If this facet does not contain a lattice
point then $\co(\Delta\cap M)$ is contained in the half-space
$\{m\in M_\RR\mid \langle m, v\rangle\ge0\}$. But then
$(\co(\Delta\cap M ))^*$ is unbounded, contradicting that $\co(
(\co(\Delta\cap M ))^*\cap N)=\Delta^*$ is a polytope.
\end{pf}

\begin{lem}\label{l:qrefprop} If $\Delta$ is a $\QQ$-reflexive polytope in $M_\RR$,
then\\
{\rm (a)} $\inte(\Delta)\cap M=\{0\}$, \\ {\rm (b)} $
\Delta^*=[\Delta^\circ]$, \\
{\rm (c)} $\inte(\Delta^*)\cap N=\{0\}$.
\end{lem}

\begin{pf}
The property {\rm (a)} follows from the fact that $\Delta^*$ is a
lattice polytope, while the property  {\rm (b)} is equation
$(\ref{e:qreflex})$ in the new notation.  Part {\rm (c)} holds
since $ \Delta^\circ$ is $\QQ$-reflexive in $N_\RR$ and applying
properties {\rm (a)}  and {\rm (b)} to $ \Delta^\circ$  we get
$\inte(\Delta^*)\cap N=\inte(\Delta^\circ)\cap N=\{0\}$.
\end{pf}

By   part  {\rm (c)} of the above lemma, we get

\begin{cor}\label{c:can} If $\Delta$ is a $\QQ$-reflexive polytope, then $\Delta^*$ is a
canonical  Fano   polytope.
\end{cor}

\section{$\QQ$-nef-partitions.}

In this section, we generalize the construction of nef-partitions
of  L. Borisov in \cite{bo} in the context of $\QQ$-reflexive
polytopes.

\begin{defn}\cite{bo} A {\it nef-partition} of a reflexive polytope
$\Delta$ is a Minkowski sum decomposition
$\Delta=\Delta_1+\dots+\Delta_r$ by lattice polytopes such that
 $0\in\Delta_i$ for all $i$.
\end{defn}

 \begin{thm}\cite{bo} Let $\Delta=\Delta_1+\dots+\Delta_r$ be a nef-partition. If
  $$\nabla_j=\{y\in N_\RR \mid\langle
 x,y\rangle\ge-\delta_{ij}\, \forall \,x\in\Delta_i,\,
 i=1,\dots,r\}$$
  for $j=1,\dots,r$, where $\delta_{ij}$ is
 the Kronecker symbol,
  then
   $\nabla=\nabla_1+\dots+\nabla_r$ is a  nef-partition. Moreover, $$\Delta_i=\{x\in M_\RR \mid\langle
 x,y\rangle\ge-\delta_{ij}\, \forall \,y\in \nabla_j,\,
 j=1,\dots,r\}$$
 for $i=1,\dots,r$.
\end{thm}

The nef-partitions $\Delta=\Delta_1+\dots+\Delta_r$ and
 $\nabla=\nabla_1+\dots+\nabla_r$ are called a {\it dual  pair of
 nef-partitions}.

\begin{rem} The name nef-partition comes from two words: nef and
partition. The {\it nef} part comes from  the property that each
summand $\Delta_i$ in the Minkowski sum
$\Delta=\Delta_1+\dots+\Delta_r$ defines a   nef (numerically
effective) divisor
$$D_{\Delta_i}=
\sum_{\rho\in\Sigma_\Delta(1)}(-\min\langle\Delta_i,v_\rho\rangle)
D_\rho= \sum_{v_\rho\in\nabla_i} D_\rho$$ on the Gorenstein Fano
toric variety $X_\Delta$, where $D_\rho$ are the torus invariant
divisors in $X_\Delta$ corresponding to the rays $\rho$ of the
normal fan $\Sigma_\Delta$ of the polytope $\Delta$, and $v_\rho$
are the primitive lattice generators of $\rho$. The {\it
partition} part corresponds to the fact that the anticanonical
divisor has its support $\bigcup_{\rho\in\Sigma_\Delta(1)}D_\rho$
  partitioned into the union of supports
$\bigcup_{v_\rho\in\nabla_i} D_\rho$ of the nef-divisors
$D_{\Delta_i}$.
\end{rem}

\begin{rem} It was an original idea of Yu.~I.~Manin (see \cite[Sect.~6.2]{bv}) to partition
the disjoint union  $\bigcup_{\rho\in\Sigma_\Delta(1)}D_\rho$ of
torus invariant divisors into a union of sets which support the
nef-divisors $D_{\Delta_i}$. L. Borisov translated this idea into
Minkowski sums and found a canonical way of creating dual
nef-partitions.
\end{rem}

 Now, we introduce a generalization of nef-partition in the context of $\QQ$-reflexive polytopes.

\begin{defn}\label{d:qnef} A {\it $\QQ$-nef-partition} of a $\QQ$-reflexive polytope
$\Delta$ is a Minkowski sum decomposition
$\Delta=\Delta_1+\cdots+\Delta_r$  into polytopes in $M_\RR$ such
that
  $0\in\Delta_i$ for all $i$, and $\co(\Delta\cap M)=\co(\Delta_1\cap M)+\cdots+\co(\Delta_r\cap
  M).$
\end{defn}

 A {\it $\QQ$-nef-partition} has the following property.

\begin{lem} Let $\Delta=\Delta_1+\cdots+\Delta_r$   be a $\QQ$-nef-partition,
and let $F$ be a facet of $\Delta$ and $F=F_1+\cdots+F_r$  be the
induced decomposition by faces $F_i$ of $\Delta_i$, for
$i=1,\dots,r$.  Then $\co(F\cap M)=\co(F_1\cap
M)+\cdots+\co(F_r\cap
  M) $.
\end{lem}

\begin{pf}  Let $F$ be a facet of  $\Delta$ with the
induced decomposition $F=F_1+\cdots+F_r$. Then the inclusion
$\co(F_1\cap M)+\cdots+\co(F_r\cap
  M)\subseteq \co(F\cap M) $ is clear.
To show the other inclusion, notice that $[F]=\co(F\cap M) $ is a
nonempty face of $\co(\Delta\cap M)$, by Lemma~\ref{l:facet}.  By
Definition~\ref{d:qnef}, we have
$[\Delta]=[\Delta_1]+\cdots+[\Delta_r]$, which induces the
Minkowski sum decomposition $[F]=G_1+\cdots+G_r$ by faces $G_i$ of
$[\Delta_i]$. Let $v$ be the vertex of $\Delta^*$ such that
$\langle F,v\rangle=\min\langle\Delta,v\rangle=-1$. We have
$\min\langle\Delta_i,v\rangle\le \min\langle[\Delta_i],v\rangle$
for all $i$ and
\begin{equation}\label{e:ineq}\min\langle\Delta,v\rangle=\sum_{i=1}^r\min\langle\Delta_i,v\rangle\le\sum_{i=1}^r
\min\langle[\Delta_i],v\rangle=\min\langle[\Delta],v\rangle.\end{equation}
Hence, $\min\langle\Delta_i,v\rangle=
\min\langle[\Delta_i],v\rangle$ for all $i$, since
$\min\langle\Delta,v\rangle= \min\langle[\Delta],v\rangle=-1$  by
 Lemma~\ref{l:facet}. Since the faces $F_i$ and $G_i$ are
 determined by the minimal value of $v$ on $\Delta_i$ and
 $[\Delta_i]$, respectively,  we conclude that $G_i\subseteq F_i$,
 whence $[F]=G_1+\cdots+G_r\subseteq \co(F_1\cap M)+\cdots+\co(F_r\cap
  M)$.
\end{pf}

\begin{defn} For a   $\QQ$-nef-partition
$ \Delta_1+\cdots+\Delta_r$ in $M_\RR$ define the polytopes
\begin{equation}\label{e:dualnef}\nabla_j=\{y\in N_\RR
\mid\langle
 x,y\rangle\ge-\delta_{ij}\, \forall \,x\in\co(\Delta_i\cap M),\,
 i=1,\dots,r\}\end{equation}
  for $j=1,\dots,r$.
\end{defn}

\begin{pr}\label{p:first} Let  $ \Delta_1+\cdots+\Delta_r$ be a
$\QQ$-nef-partition in $M_\RR$,  then
$$(\Delta_1+\cdots+\Delta_r)^*= \co(\nabla_1\cap
N,\dots,\nabla_r\cap N),$$   where $\nabla_1,\dots,\nabla_r$ are
defined by $(\ref{e:dualnef})$.
\end{pr}

\begin{pf} Let $v$ be a vertex of $\Delta^*$, where
$\Delta=\Delta_1+\cdots+\Delta_r$ is a $\QQ$-nef-partition. Then
by $(\ref{e:ineq})$ and Lemma~\ref{l:facet}, we have $\sum_{i=1}^r
\min\langle[\Delta_i],v\rangle=-1$, whence the integer
$\min\langle[\Delta_j],v\rangle=-1$ for some $j$ and
$\min\langle[\Delta_i],v\rangle=0$ for $i\ne j$ since
$\min\langle[\Delta_i],v\rangle\le0$ by $0\in[\Delta_i]$ for all
$i$. Hence, every vertex $v$ of $\Delta^*$ is contained in some
$\nabla_j\cap N$, and $\Delta^*\subseteq  \co(\nabla_1\cap
N,\dots,\nabla_r\cap N) $.

To show the opposite inclusion, let $y\in \nabla_j\cap N$ for some
$j$. Then $\min\langle[\Delta],y\rangle= \sum_{i=1}^r
\min\langle[\Delta_i],y\rangle\ge\sum_{i=1}^r\delta_{ij}=-1$,
whence $y\in [\Delta]^*=\Delta^\circ$. Since $y$ is a lattice
point, we get $y\in [\Delta^\circ]=\Delta^*$ by part {\rm (b)} of
Lemma~\ref{l:qrefprop}.
\end{pf}

\begin{pr}\label{p:2nd} Let  $\Delta_1+\cdots+\Delta_r$ be a
$\QQ$-nef-partition in $M_\RR$,  then $$
(\nabla_1+\cdots+\nabla_r)^*=\co(\Delta_1\cap M,\dots,\Delta_r\cap
M),$$    where $\nabla_1,\dots,\nabla_r$ are defined by
$(\ref{e:dualnef})$.
\end{pr}

\begin{pf}
 One inclusion $\co([\Delta_1],\dots,[\Delta_r])\subseteq(\nabla_1+\cdots+\nabla_r)^*$ holds since $\langle
[\Delta_i],\nabla_j\rangle\ge-\delta_{ij}$ for all $i,j$ by
$(\ref{e:dualnef})$.

The opposite inclusion holds because  $[\Delta_i]$ are lattice
polytopes with $0\in[\Delta_i]$ and $0$ is the only interior
lattice point in $[\Delta_1]+\cdots+[\Delta_r]$ by
Definition~\ref{d:qnef} and part $ \rm( a)$ of
Lemma~\ref{l:qrefprop}
\end{pf}

\begin{defn}  Let $P_1,\dots,P_r$ be polytopes in
$M_\RR$. Consider the lattice $\bar M=M\oplus\ZZ^r$, where
$\{e_1,\dots,e_r\}$ is the standard basis of $\ZZ^r$. The cone
$$  {\cal C}_{P_1,\dots,P_r}:=\RR_{\ge0}\cdot{\rm Conv}(
   P_1+e_1,\dots,P_r +e_r) $$ is called the
  {\it Cayley cone} associated to the $r$-tuple of polytopes
  $P_1,\dots,P_r$.
\end{defn}

\begin{lem} \cite[Lem.~1.6]{m3}   Let $P_1,\dots,P_r$ be   polytopes in
$M_\RR$ such that $P=P_1+\cdots+P_r$ is $d$-dimensional and
$0\in\inte(P)$ Then the dual of the Cayley cone associated to
$P_1,\dots,P_r$ is
$${\cal C}_{P_1,\dots,P_r}^\ve=
\RR_{\ge0}\cdot\co\biggl(\biggl\{x-\sum_{i=1}^r \min\langle
 \Delta_i,x\rangle
 e^*_i\mid x\in \cv(P^*)\biggr\},\{e^*_1, \dots,e^*_r\}\biggr),$$
 where $\cv(P^*)$ is the set of vertices of $P^*$ and $\{e^*_1, \dots,e^*_r\}$ is the   basis of $\ZZ^{r}\subset N\oplus\ZZ^{r}$ dual to
$\{e_1, \dots,e_r\}$.
\end{lem}

Using this lemma from Propositions~\ref{p:first} and \ref{p:2nd}
we get

\begin{pr}\label{p:caydual} Let  $ \Delta_1+\cdots+\Delta_r$ be a
$\QQ$-nef-partition in $M_\RR$,  then
$${\cal C}_{\Delta_1,\dots,\Delta_r}^\ve={\cal C}_{[\nabla_1],\dots,[\nabla_r]},
\quad {\cal C}_{\nabla_1,\dots,\nabla_r}^\ve={\cal
C}_{[\Delta_1],\dots,[\Delta_r]},$$  where
$\nabla_1,\dots,\nabla_r$ are defined by $(\ref{e:dualnef})$.
\end{pr}

Applying a projection technique from \cite{bn}  to the Cayley
cones in the last proposition we get the following one.

\begin{pr}\label{p:3rd} Let  $ \Delta_1+\cdots+\Delta_r$ be a
$\QQ$-nef-partition in $M_\RR$,  then
$$(\co(\Delta_1,\dots,\Delta_r))^*=\co(\nabla_1\cap
N)+\cdots+\co(\nabla_r\cap N),$$
$$(\co(\nabla_1,\dots,\nabla_r))^*=\co(\Delta_1\cap
M)+\cdots+\co(\Delta_r\cap M),$$
   where $\nabla_1,\dots,\nabla_r$ are defined by
$(\ref{e:dualnef})$.
\end{pr}

\begin{pr} Let  $ \Delta_1+\cdots+\Delta_r$ be a
$\QQ$-nef-partition in $M_\RR$, and let $\nabla_j$ be defined by
$(\ref{e:dualnef})$.   Then $\co(\nabla_1,\dots,\nabla_r)$ is a
$\QQ$-reflexive polytope  and
$$ (\Delta_1+\cdots+\Delta_r)^\circ=\co(\nabla_1,\dots,\nabla_r).$$
\end{pr}

\begin{pf} By Definitions~\ref{d:circ}, \ref{d:qnef},  and  Proposition~\ref{p:3rd}  we have
$$(\Delta_1 +\cdots+ \Delta_r)^\circ=[\Delta_1 +\cdots+
\Delta_r]^*=([\Delta_1]+\cdots+[\Delta_r])^*=
\co(\nabla_1,\dots,\nabla_r).$$\end{pf}

\begin{pr}\label{p:2ndqp}  Let  $ \Delta_1+\cdots+\Delta_r$ be a
$\QQ$-nef-partition in $M_\RR$, and let $\nabla_1\dots,\nabla_r$
be defined by $(\ref{e:dualnef})$.   Then
$\nabla_1+\cdots+\nabla_r$ and $\co(\Delta_1,\dots,\Delta_r)$  are
$\QQ$-reflexive polytopes  and
$$(\nabla_1+\cdots+\nabla_r)^\circ=\co(\Delta_1,\dots,\Delta_r).$$
\end{pr}

\begin{pf} First, we claim
\begin{equation}\label{e:equallat} \co([\Delta_1],\dots,[\Delta_r])=[\co(\Delta_1,\dots,\Delta_r)].\end{equation}
Indeed, take $x\in \co(\Delta_1,\dots,\Delta_r)\cap M$, then
$x\in(\sum_{i=1}^r[\nabla_i])^*$ by Proposition~\ref{p:3rd}, and
we have $$-1\le\min\langle x,\sum_{i=1}^r[\nabla_i]
\rangle=\sum_{i=1}^r\min\langle x,[\nabla_i] \rangle\le0,$$ since
$0\in [\nabla_i]$ for all $i$. From here we get two cases: either
all integers $\min\langle x,[\nabla_i] \rangle=0$, for
$i=1,\dots,r$, or there is $j$ such that $\min\langle x,[\nabla_i]
\rangle=-\delta_{ij}$ for $i=1,\dots,r$. This implies that either
$x=0$ or $x\in\Delta_j$.  Thus, we showed that
$(\ref{e:equallat})$ holds.

Next, we show that $\nabla_1 +\cdots+ \nabla_r$ is
$\QQ$-reflexive. Using Propositions~\ref{p:2nd}, \ref{p:3rd}, we
have
\begin{align*}\co([\Delta_1],\dots,[\Delta_r])=&(\nabla_1 +\cdots+
\nabla_r)^*\subseteq  [\nabla_1 +\cdots+
\nabla_r]^*\subseteq\\&\subseteq ([\nabla_1] +\cdots+
[\nabla_r])^*=\co(\Delta_1,\dots,\Delta_r).\end{align*} Applying
$(\ref{e:equallat})$ to this, we get $$[[\nabla_1 +\cdots+
\nabla_r]^*]=\co([\Delta_1],\dots,[\Delta_r])=(\nabla_1 +\cdots+
\nabla_r)^*,$$  showing that the polytope $\nabla_1 +\cdots+
\nabla_r$ is $\QQ$-reflexive. Then, by Definition~\ref{d:circ} and
the properties of $\QQ$-reflexive polytopes, the dual
$\QQ$-reflexive polytope is $(\nabla_1 +\cdots+
\nabla_r)^\circ=\co(\Delta_1,\dots,\Delta_r)$.
\end{pf}

Finally, we establish the existence of the dual
$\QQ$-nef-partition:

\begin{thm}\label{t:qnefex} Let  $\Delta_1+\cdots+\Delta_r$ be a
$\QQ$-nef-partition, then $ \nabla_1+\cdots+\nabla_r$ is   a
$\QQ$-nef-partition,  where $\nabla_1,\dots,\nabla_r$ are defined
by $(\ref{e:dualnef})$. Moreover, $$\Delta_i=\{x\in M_\RR
\mid\langle
 x,y\rangle\ge-\delta_{ij}\, \forall \,y\in\co(\nabla_j\cap N),\,
 j=1,\dots,r\}.$$
\end{thm}

\begin{pf} Proposition~\ref{p:2ndqp} gives $\QQ$-reflexivity of   $ \nabla_1+\cdots+\nabla_r$. To show that  $
\nabla_1+\cdots+\nabla_r$ is a $\QQ$-nef-partition notice
$$(\co(\Delta_1,\dots,\Delta_r))^*=[\nabla_1]+\cdots+[\nabla_r]\subseteq[\nabla_1+\cdots+\nabla_r]$$
by Proposition~\ref{p:3rd}. Applying part (b) of
Lemma~\ref{l:qrefprop}, we see that
$[\nabla_1+\cdots+\nabla_r]=(\co(\Delta_1,\dots,\Delta_r))^*$
since
$\co(\Delta_1,\dots,\Delta_r)=(\nabla_1+\cdots+\nabla_r)^\circ$ by
Proposition~\ref{p:2ndqp}. From the above inclusions we get  the
required equality
$[\nabla_1+\cdots+\nabla_r]=[\nabla_1]+\cdots+[\nabla_r]$ in the
definition of a $\QQ$-nef-partition. The last part of this theorem
follows by Proposition~\ref{p:first} since $0$ is the only
interior lattice point in $[\nabla_1]+\cdots+[\nabla_r]$.
\end{pf}

The Minkowski sums $\Delta_1+\dots+\Delta_r$ and
 $\nabla_1+\dots+\nabla_r$ in Theorem~\ref{t:qnefex} will be called a {\it dual  pair of
$\QQ$-nef-partitions}.

\begin{defn} A $\QQ$-nef-partition $\Delta_1+\dots+\Delta_r$ in $M_\RR$ is
called {\it proper} if $ \Delta_i\ne 0 $  for all $1\le i\le r$.
\end{defn}

\begin{cor} Let  $\Delta_1+\dots+\Delta_r\subset M_\RR$ and
 $\nabla_1+\dots+\nabla_r\subset N_\RR$ be a dual  pair of
$\QQ$-nef-partitions. Then, for $i=1,\dots,r$, one has $ \Delta_i
= 0 $ if and only if $ \nabla_i   = 0 $.
\end{cor}

\begin{pf} If $
\Delta_i = 0 $, then $\nabla_i=0 $ by Definition~\ref{d:qnef},
since $[\Delta_1]+\cdots+[\Delta_r]$ spans $M_\RR$. The opposite
implication follows from Theorem~\ref{t:qnefex}.
\end{pf}

\begin{cor} If $\Delta_1+\dots+\Delta_r\subset M_\RR$ is a proper
$\QQ$-nef-partition, then its dual $\QQ$-nef-partition
$\nabla_1+\dots+\nabla_r$ is proper.
\end{cor}

\begin{cor} If $\Delta_1+\dots+\Delta_r\subset M_\RR$ is a proper
$\QQ$-nef-partition, then  $ \Delta_i\cap M\ne\{0\}$ for all $1\le
i\le r$.
\end{cor}

\section{Almost reflexive Gorenstein cones.}

In this section, we generalize the notion of reflexive Gorenstein
cones.

\begin{defn}\cite{bb} Let $\bar M$ and $\bar N$ be a pair of   dual lattices of rank $\bar d$.
    A $\bar
d$-dimensional polyhedral cone
 $\sigma$ with a vertex
at $0\in \bar M$ is called {\em Gorenstein}, if  it is generated
by finitely many lattice points contained in the affine hyperplane
$\{x\in \bar M\mid
 \langle x,h_\sigma\rangle=1\}$ for   $h_{\sigma}\in \bar N$.
 The unique lattice point
$h_\sigma$
 is called the {\it height} (or {\it degree}) vector of the Gorenstein cone $\sigma$.
  A Gorenstein cone $\sigma$ is called   {\it reflexive}
 if both $\sigma$ and its dual
$$\sigma^\ve=\{y\in \bar N_\RR\mid \langle x, y\rangle\ge0
\,\forall\,  x\in\sigma \}$$ are Gorenstein cones.   In this case,
they both have uniquely determined $h_\sigma\in  \bar N$ and
$h_{\sigma^\ve}\in \bar M$, which take value 1 at the primitive
lattice generators of the respective cones.  The positive integer
$r=\langle h_{\sigma^\ve},h_\sigma\rangle$ is called  the {\it
index} of the reflexive Gorenstein cones $\sigma$ and
$\sigma^\ve$.
\end{defn}

As in \cite{bn}, denote  $\sigma_{(i)}:=\{x\in\sigma\mid\langle
x,h_\sigma\rangle=i\}$, for $i\in\NN$.  The basic relationship
between reflexive polytopes and reflexive Gorenstein cones is
provided by the following:

\begin{pr}\cite[Pr.~2.11]{bb}\label{p:reflg} Let $\sigma$ be a Gorenstein cone. Then
$\sigma$ is a reflexive Gorenstein cone of index $r$ if and only
if the polytope $\sigma_{(r)}-h_{\sigma^\ve}$ is a reflexive
polytope with respect to the lattice $\bar M\cap
h_\sigma^\perp=\{x\in \bar M\mid \langle x,h_\sigma\rangle=0\}$.
\end{pr}

Generalizing the notion of reflexive Gorenstein cones we
introduce:

\begin{defn}  A  Gorenstein cone
 $\sigma$  in $\bar M_\RR$ is called {\it almost reflexive}, if
 there is $r\in\NN$ such that  $\sigma_{(r)}$ has a unique lattice
 point $h$ in its relative interior and $\sigma_{(r)}-h$ is an
 almost reflexive polytope with respect to the lattice $\bar M\cap
h_\sigma^\perp$. We will denote $h$ by $h_{\sigma^\ve}$. The
positive integer $r$ will be called the {\it  index} of the almost
reflexive Gorenstein cone  $\sigma$.
\end{defn}

\begin{lem} Reflexive Gorenstein cones are almost reflexive.
\end{lem}

\begin{pf} This follows from Proposition~\ref{p:reflg} and
Lemma~\ref{qrefaref}.
\end{pf}

\begin{defn} For an almost reflexive Gorenstein cone  $\sigma$  in $\bar M_\RR$ define
$$\sigma^\ve_{(i)}=\{y\in\sigma^\ve\mid\langle
 h_{\sigma^\ve},y\rangle=i\}, \,\text{ for }\, i\in\NN.$$
 Denote $[\sigma^\ve]:=\RR_{\ge0}[\sigma^\ve_{(1)}]=\RR_{\ge0}\co(\sigma^\ve_{(1)}\cap \bar M)$
\end{defn}

\begin{lem}\label{l:dual1} Let $\Delta$ be a polytope in $M_\RR$ with
$0\in\inte(\Delta)$, and let
 $\sigma_\Delta=\RR_{\ge0}(\Delta,1)\subset \bar
M_\RR=M_\RR\oplus\RR$.
 Then $$\sigma_\Delta^\ve=\sigma_{\Delta^*}=\RR_{\ge0}(\Delta^*,1)\subset \bar
N_\RR=N_\RR\oplus\RR.$$
\end{lem}

\begin{cor}
A Gorenstein cone $\sigma$ in $\bar M_\RR$ is  almost reflexive of
index 1 if and only if the polytope $\sigma^\ve_{(1)}-h_\sigma$ is
$\QQ$-reflexive with respect to the lattice $\bar N\cap
h_{\sigma^\ve}^\perp$.
\end{cor}

\begin{pf} Combine    Lemmas~\ref{qrefaref} and \ref{l:dual1}.
\end{pf}

\begin{cor}  If
 $\sigma$  in $\bar M_\RR$ is  an almost reflexive    Gorenstein
 cone of index 1, then $[\sigma^\ve]$ is  an almost reflexive    Gorenstein
 cone of index 1.
\end{cor}

\begin{pr}  If
 $\sigma$  in $\bar M_\RR$ is  an almost reflexive    Gorenstein
 cone of index $r$, then $[\sigma^\ve]$ is  an almost reflexive    Gorenstein
 cone of index $r$.
\end{pr}

\begin{pf} Use the techniques in the proof of \cite[Pr.~2.11]{bb}.
\end{pf}

Almost reflexive    Gorenstein
 cones have the following property.

\begin{lem} Let $\sigma\subset\bar M_\RR$ be  an almost reflexive    Gorenstein
 cone. Then
 $[[\sigma^\ve]^\ve]=\sigma$.
\end{lem}

\begin{defn} For an almost reflexive Gorenstein cone $\sigma$,
denote $\sigma^\bullet:=[\sigma^\ve]$.
\end{defn}

\begin{cor} The map $\sigma\mapsto \sigma^\bullet$ is  an involution on the set  of almost reflexive
    Gorenstein cones: $(\sigma^\bullet)^\bullet=\sigma$.
\end{cor}

 Cayley cones corresponding to a   dual pair of
$\QQ$-nef-partitions are related to almost reflexive Gorenstein
cones as follows:

\begin{pr} Let  $\Delta_1+\dots+\Delta_r\subset M_\RR$ and
 $\nabla_1+\dots+\nabla_r\subset N_\RR$ be a dual  pair of
$\QQ$-nef-partitions. Then the Cayley cones ${\cal
C}_{[\Delta_1],\dots,[\Delta_r]}$ and ${\cal
C}_{[\nabla_1],\dots,[\nabla_r]}$ is a dual pair of almost
reflexive Gorenstein cones: $${\cal
C}_{[\Delta_1],\dots,[\Delta_r]}^\bullet=[{\cal
C}_{[\Delta_1],\dots,[\Delta_r]}^\ve]={\cal
C}_{[\nabla_1],\dots,[\nabla_r]}.$$
\end{pr}

\begin{pf} This follows directly from Proposition~\ref{p:caydual}
since the height vectors of the Cayley cones ${\cal
C}_{[\Delta_1],\dots,[\Delta_r]}$ and ${\cal
C}_{[\nabla_1],\dots,[\nabla_r]}$ are $e_1^*+\cdots+e_r^*$  and
$e_1 +\cdots+e_r $, respectively.
\end{pf}

\section{Basic    toric geometry.}\label{s:basic}

This section will review some basics of toric geometry.

 Let $\xs$ be a $d$-dimensional toric variety associated
with a finite rational polyhedral fan $\Sigma$ in $N_\RR$. Denote
by $\Sigma(1)$ the finite set of the 1-dimensional cones $\rho$ in
$\Sigma$, which correspond to the torus invariant divisors
$D_\rho$ in $\xs$.    By \cite{c}, every toric variety can be
described as a categorical quotient of a Zariski open subset of an
affine space by a subgroup of a torus.  Consider the polynomial
ring $S(\Sigma):=\CC[x_\rho: \rho\in\Sigma(1)]$, called the {\it
homogeneous coordinate ring} of the toric variety $\xs$, and the
corresponding affine space $\CC^{\Sigma(1)}=\spe(\CC[x_\rho:
\rho\in\Sigma(1)])$. The ideal  $B =\langle
\prod_{\rho\not\subseteq\sigma} x_\rho: \sigma\in\Sigma\rangle$
 in $S$ is called the irrelevant ideal. This ideal determines a Zariski closed set ${\bf
V}(B )$ in $\CC^{\Sigma(1)}$, which
 is
invariant under the diagonal group action of the subgroup
$$G=\biggl\{(\mu_\rho) \in(\CC^*)^{\Sigma(1)}\mid\prod_{\rho\in\Sigma(1)} \mu_\rho^{\langle
u,v_\rho \rangle}=1\, \forall\, u\in M\biggr\}$$ of the torus
$(\CC^*)^{\Sigma(1)}$ on the affine space $\CC^{\Sigma(1)}$, where
$v_\rho$ denotes the primitive lattice generator of the
1-dimensional cone $\rho$. The toric variety $\xs$ is isomorphic
to the categorical quotient $(\CC^{\Sigma(1)}\setminus{\bf
V}(B))/G$, which is induced by a
  toric morphism $\pi:\CC^{\Sigma(1)}\setminus
Z(\Sigma)\rightarrow\xs $, constant on  $G(\Sigma)$-orbits (see
 \cite[Thm.~5.1.10]{torvar}).

The coordinate ring $S(\Sigma)$ is graded by the the Chow group
$$A_{d-1}(\xs)\simeq \Hom(G,\CC^*),$$
 and  $\deg(\prod_{\rho\in\Sigma(1) }
x_\rho^{b_\rho})= [\sum_{\rho\in\Sigma(1)} b_\rho D_\rho]\in
A_{d-1}(\xs)$. For a torus invariant Weil divisor
$D=\sum_{\rho\in\Sigma(1)} b_\rho D_\rho$, there is a one-to-one
correspondence between the monomials of $\CC[x_\rho:
\rho\in\Sigma(1)]$ in the degree $[\sum_{\rho\in\Sigma(1)} b_\rho
D_\rho]\in A_{d-1}(\xs)$ and the lattice points inside the
polytope
$$ \Delta_D=\{m\in M_\RR\mid\langle m,v_\rho\rangle\ge-b_\rho\,\forall\,
\rho\in\Sigma(1)\}
$$ by associating  to $m\in\Delta_D$ the monomial  $ \prod_{\rho\in\Sigma(1) }
x_\rho^{b_\rho+\langle m,v_\rho\rangle}$. If we denote the
homogeneous degree of $S(\Sigma)$ corresponding to $\beta=[D]\in
A_{d-1}(\xs)$ by $S(\Sigma)_\beta$, then by  \cite[Prop.~1.1]{c},
we also have a natural isomorphism
$$H^0(\xs,{\cal O}_{\xs}(D))\simeq S(\Sigma)_\beta.$$
In particular, every hypersurface in $\xs$ of degree
$\beta=\sum_{\rho\in\Sigma(1)} b_\rho D_\rho$ corresponds to a
polynomial $$\sum_{m\in\Delta_D\cap M} a_m \prod_{\rho\in\Sigma(1)
} x_\rho^{b_\rho+\langle m,v_\rho\rangle}$$ with the coefficients
$a_m\in \CC$.
 By \cite[Prop.~5.2.8]{torvar}, all closed
subvarieties of $\xs$ correspond to homogeneous ideals $I\subseteq
S(\Sigma)$, and \cite[Thm.~1.2]{m2} shows that a closed subvariety
in a toric variety $\xs$ can be viewed as a categorical quotient
  as well. A {\it complete intersection} in the toric variety $\xs$ (in homogeneous coordinates)
  is
  a
closed  subvariety ${\bf V}(I)\subset\xs$ corresponding to a
radical homogeneous ideal $I\subseteq S(\Sigma)$ generated
  by a  regular sequence of homogeneous polynomials $f_1,\dots,f_k\in S(\Sigma)$ such that $k=\dim
  \xs-\dim {\bf V}(I)$ (see \cite[Sect.~1]{m2}).

Every rational polytope $\Delta$ in $M_\RR$ determines the Weil
$\QQ$-divisor
$$D_\Delta=\sum_{\rho\in\Sigma(1)} (-\min\langle
\Delta,v_\rho\rangle) D_\rho\in{\rm WDiv}(\xs)\otimes_\ZZ\QQ$$ on
$\xs$.

\begin{defn} Let $X$ be a complete variety. A $\QQ$-Cartier divisor $D\in{\rm Div}(X)\otimes_\ZZ\QQ$ on
$X$ is called ${\it nef}$ (numerically effective) if $D\cdot
C\ge0$ for all irreducible curves $C\subset X$. We will call such
divisors {\it $\QQ$-nef}.
\end{defn}

\begin{lem} Let $\xs$ be a compact toric variety. Then  the divisor $D_\Delta$ is
$\QQ$-nef if and only if its support function
$\psi_\Delta=-\min\langle\Delta,\underline{\,\,\,}\,\rangle$ is
 convex piecewise linear with respect to the fan
$\Sigma$.
\end{lem}

\begin{pf} By \cite[p.~68]{f}, we know that ${\cal
O}_\xs(nD_\Delta)$ is generated by global sections for some
sufficiently large $n\in\NN$ if and only if $\psi_\Delta$ is
 convex piecewise linear on
$\Sigma$. On the other hand, we showed in \cite[Thm.~1.6]{m0}
that, for a compact toric variety $\xs$, the invertible sheaf
${\cal O}_\xs(D)$ is generated by global sections if and only if
$D$ is nef.
\end{pf}

\begin{lem}\cite[Lem~2.1]{m2}\label{l:sumnef} Let $\xs$ be a compact toric variety   associated to a fan $\Sigma$ in
$N_\RR$. Suppose ${\Delta_1}$ and $ {\Delta_2}$ are rational
polytopes in $M_\RR$ then $D_{\Delta_1+\Delta_2}$ is a $\QQ$-nef
divisor on $\xs$ iff $D_{\Delta_1}$ and $D_{\Delta_2}$ are
$\QQ$-nef on $\xs$.
\end{lem}

Every rational polytope $\Delta$ in $M_\RR$ corresponds to a
projective toric variety $\xd:=\xsd$, whose fan $\Sigma_\Delta$
(called the {\it normal fan of $\Delta$})  is the collection of
cones
$$\sigma_F=\{y\in N_\RR \mid \langle
x,y\rangle\le\min\langle\Delta,y\rangle\,\forall\,x\in F\}.$$

The support function
$\psi_\Delta=-\min\langle\Delta,\underline{\,\,\,}\,\rangle$ is
strictly convex piecewise linear with respect to the fan
$\Sigma_\Delta$. In this case, the divisor $D_\Delta$ is ample,
and, in particular, $\QQ$-nef. From Lemma~\ref{l:sumnef}, we get

\begin{cor} Let $\xd$ be a Fano toric variety, and suppose
$\Delta=\Delta_1+\cdots+\Delta_r$ is a Minkowski sum decomposition
by rational polytopes. Then the divisors $D_{\Delta_i}$ are
$\QQ$-nef  on $\xd$ for all $1\le i\le r$.
\end{cor}

There is an alternative way to describe projective toric varieties
using the Proj functor, which is simple but less useful in the
context of complete intersections. Consider the   cone $$K=\{(
t\Delta,t)\mid t\in\RR_{\ge0}\}\subset M_\RR\oplus\RR.$$
 The projective toric variety
 $\xd=\xsd$
can be represented as $ {\rm Proj}({\Bbb C}[K\cap(M\oplus\ZZ)]).$
 Moreover, if  $\beta\in A_{d-1}(\xd)$ is
the class of the  ample divisor
$D_\Delta=\sum_{\rho\in\Sigma_\Delta(1)} b_\rho D_\rho$, then
there is a natural isomorphism of graded rings
$$
{\Bbb C}[K\cap(M\oplus\ZZ)] \simeq\bigoplus_{i=0}^\infty
S(\sd)_{i\beta},
$$
sending $\chi^{(m,i)}\in\CC[K\cap(M\oplus\ZZ)]_i$ to
$\prod_{\rho\in\sd(1)} x_\rho^{i b_\rho+\langle m,v_\rho\rangle}$.
In particular, a hypersurface given by a polynomial  in
homogeneous coordinates
$$\sum_{m\in\Delta\cap M}a_m \prod_{\rho\in\sd(1)} x_\rho^{i
b_\rho+\langle m,v_\rho\rangle}= 0
$$ corresponds to $\sum_{m\in\Delta\cap M}a_m \chi^{(m,i)}=0$.

\section{Mirror Symmetry Construction.}

In this section, we propose a generalization of the
Batyrev-Borisov Mirror Symmetry constructions.

A proper $\QQ$-nef-partition $\Delta=\Delta_1+\cdots+\Delta_r$
with $r<d$ defines a $\QQ$-nef Calabi-Yau complete intersection
$Y_{\Delta_1,\dots,\Delta_r}$ in the Fano toric variety $X_\Delta$
given by the equations:
$$
 \Biggl(\sum_{m\in  \Delta_i\cap M  }a_{i,m}
\prod_{v_\rho\in\cv(\Delta^*)}x_\rho^{ \langle m,v_\rho\rangle
}\Biggr)\prod_{
 v_\rho \in [\nabla_i] } x_\rho=0,\quad i=1,\dots,r,$$
where $x_\rho$ are the  homogeneous coordinates of the toric
variety $X_\Delta$ corresponding   to the vertices $v_\rho$ of the
polytope $\Delta^*$.

Following the Batyrev-Borisov Mirror Symmetry construction,
  we naturally expect that Calabi-Yau complete intersections
corresponding to a dual pair of $\QQ$-nef-partitions pass the
topological mirror symmetry test:

\begin{conj}  Let $Y_{\Delta_1,\dots,\Delta_r}\subset X_\Delta$ and $Y_{\nabla_1,\dots,\nabla_r}\subset X_\nabla$ be a
  pair of  generic Calabi-Yau complete intersections in $d$-dimensional
Fano toric varieties corresponding to a dual pair of
$\QQ$-nef-partitions $\Delta=\Delta_1+\cdots+\Delta_r$ and
$\nabla=\nabla_1+\cdots+\nabla_r$. Then
$$h^{p,q}_{\rm st}(Y_{\Delta_1,\dots,\Delta_r})=h^{d-r-p,q}_{\rm
st}(Y_{\nabla_1,\dots,\nabla_r}),\quad 0\le p,q\le d-r.$$
\end{conj}

Assuming that this conjecture holds, by taking maximal projective
crepant partial resolutions  we obtain the mirror pair of minimal
Calabi-Yau complete intersections $\widehat
Y_{\Delta_1,\dots,\Delta_k}$, $\widehat
Y_{\nabla_1,\dots,\nabla_k}$.

For an almost reflexive Gorenstein cone $\sigma$ in $\bar M_\RR$,
we have the Fano toric variety $X_\sigma=\pro(\CC[\sigma^\ve\cap
\bar N])$, whose fan consists of cones generated by the faces of
the almost reflexive polytope $\sigma_{(r)}-h_{\sigma^\ve}$ in
$\bar M_\RR\cap h_\sigma^\perp$. A {\it generalized Calabi-Yau
manifold} is defined as the ample $\QQ$-Cartier hypersurface
$Y_\sigma\subset X_\sigma$ given by the equation
$$\sum_{n\in\sigma^\bullet_{(1)}\cap \bar N}a_n\chi^n=0$$ with
generic $a_n\in\CC$, where $\chi^n$ are the elements in the graded
semigroup ring $\CC[\sigma^\ve\cap  \bar N]$ corresponding to the
lattice points $n\in \sigma^\ve\cap  \bar N$.

\begin{conj}  The involution $\sigma\mapsto \sigma^\bullet$ on the set of almost reflexive
    Gorenstein cones
corresponds to    the mirror involution of  $N=2$ super conformal
field theories associated to the generalized Calabi-Yau manifolds
$Y_\sigma$ and $Y_{\sigma^\bullet}$.
\end{conj}

In the case, when  a $\QQ$-nef Calabi-Yau complete intersection
$Y_{\Delta_1,\dots,\Delta_r}$  does not have the property that
$0\in\Delta_i$ for all $1\le i\le r$ (i.e, the Minkowski sum
$\Delta_1+\cdots+\Delta_r$
 is not a
$\QQ$-nef-partition), one can still associate to it  the mirror in
the form of the generalized Calabi-Yau manifold corresponding to
the dual of the Cayley cone ${\cal C}_{\Delta_1,\dots,\Delta_r}$.

\end{document}